\documentclass[11pt]{amsart}
\usepackage{amsmath}
\usepackage{latexsym}
\usepackage{amssymb}
\newtheorem{theorem}{Theorem}[section]

\newtheorem{lemma}[theorem]{Lemma}
\newtheorem{assertion}[theorem]{Assertion}

\author{Ron Aharoni}
\address{Department of Mathematics\\Technion, Haifa\\ Israel 32000}
\email[Ron Aharoni]{ra@tx.technion.ac.il}
\thanks{The research of the first author was 
supported by the fund for the promotion of research at the Technion}

\author{Eli Berger}
\email[Eli Berger]{seli@t2.technion.ac.il}

\title{The number of edges in critical strongly connected graphs }
\begin{document}
\maketitle

\begin{abstract}

We prove that the maximal number of  directed edges in a vertex-critical
strongly connected simple digraph on $n$ vertices
is $\binom{n}{2} - n +4$.
\end {abstract}

\section{Introduction}

A directed graph (or {\em digraph}) without loops or multiple edges
is called {\em strongly connected}
if each vertex in it is reachable from every other vertex. It is called
(vertex)
{\em critical strongly connected} 
if, in addition to being strongly connected, 
it has the property that the removal of any
vertex from it results in a non-strongly connected graph. 
We denote by $M(n)$ 
the maximal number of 
edges in a critical strongly connected digraph on $n$ vertices. 
Schwarz \cite{schwarz} conjectured (and proved for $n \leq 5$)
that $M(n) \leq \binom{n}{2}$. This conjecture was proved 
by London in \cite{london}. In this paper we  determine the precise
number of $M(n)$, showing that it is 
$\binom{n}{2} - n +4$. (The corresponding number for 
edge-critical strongly connected graphs is $2n-2$, see
e.g. \cite{brualdiryser}, pp 65-66.)

Here is some notation we shall use. Given a digraph $D$ we
denote by $V(D)$ the set of its vertices, and by $E(D)$
the set of edges. 
Throughout the paper the notation $n$ will be reserved for 
the number of vertices in the digraph named $D$.
For a vertex $v$ of $D$ we write $E^+_D(v)$ for the 
set of vertices $u$ for which $(v,u) \in E(D)$
and $E^-_D(v)$ for the 
set of vertices $u$ for which $(u,v) \in E(D)$.
We write $d_D(v)$ for the degree of  $v$, namely $|E^+_D(v)|+ |E^-_D(v)|$.
For a subset $A$ of $V(D)$ we write $D -A$ for
the graph obtained from $D$ by
removing all vertices in $A$, together with all edges incident with them.
If $A$ consists of a single vertex $a$, we write $D-a$ for $D - \{a\}$.
By $D/A$ we denote the digraph 
obtained from $D$ by contracting $A$, namely replacing all vertices
of $A$ by a single vertex $a$, and defining
$E^+_{D/A}(a) = \bigcup \{E^+_D(v): v \in A\} \setminus A$ 
and $E^-_{D/A}(a) = \bigcup \{E^-_D(v): v \in A\} \setminus A$.

\section{The number of edges in vertex-critical graphs}

\begin{theorem}
\label{main}
For $n \geq 4$
$$M(n) = \binom{n}{2} - n +4$$
\end{theorem}

A vertex-critical graph with  
$\binom{n}{2} - n +4$ edges is the following.
Take a directed cycle $(v_1, v_2, \ldots ,v_n)$, and add the directed edges
$(v_i, v_j), ~3\leq j < i \leq n$ and the edge $(v_2, v_1)$. Thus, what 
remains to be proved 
is that in a vertex-critical graph the
number of edges does not exceed $\binom{n}{2} - n +4$.
\\

The proof will be based on two lemmas.

\begin{lemma}
\label{theorem1}
Let $D$ be a strongly connected digraph and $v \in V(D)$ a vertex 
satisfying $d(v) \geq n$. Then there exists a vertex $z \in V(D)
\setminus \{v\}$ 
such that  $D - z$ is strongly connected.
\end{lemma} 

{\em Proof}~
The proof is by induction on $n$.
For $n<2$ the lemma is vaccuously true, since its conditions are impossible
to fulfil.
For $n=2$  take $z$ to be the vertex of the graph different from $v$. 
Let now $n>2$ and suppose that the lemma is true for all graphs
with fewer than $n$ vertices. Let $v$ be as in the lemma.
There exists then a vertex $u$ such that between $u$ and $v$ there is a
double-arc (that is, two oppositely directed eges).
Let $C = D/\{u,v\}$, and name $w$ the vertex  of $C$
replacing the shrunk pair $\{u,v\}$. By a negation hypothesis, we may assume
that $D-u$ is not strongly connected. We claim then that
$d_C(w) \geq n-1$. This will prove the lemma, since by the induction
hypothesis it will follow that $C$ has a vertex $z$ different from $w$ whose
removal leaves $C$ strongly connected. But then, clearly, also $D-z$ is
strongly connected.

To prove the claim note, first, that
$d_C(w) \geq n-2$. This follows from the fact 
that each edge in $D$ incident with $v$ and different from the two
edges joining $v$ with $u$, 
 has its copy in $C$.
Since by our assumption $D-u$ is not strongly connected, there are two edges
in $D$, say $(x,u)$ and $(u,y)$, such that $y$ is not reachable from $x$ in
$D-u$.
If $x=v$ then 
the edge $(w,y)$ is an edge in $C$ not having a 
copy $(v,y)$ in $D$, and thus can be added to the $n-2$ edges incident with
$v$ counted above, and thus $d_C(w) \geq n-1$, as desired.  
Similarly,
if $y=v$ then the edge $(x,w)$ shows that $d_C(w) \geq n-1$.
If, on the other hand, 
$x \neq v$
and $y \neq v$, then 
one of $(x,w)$ or $(w,y)$ is an edge in $C$ not counted above. 
{\hfill $\Box$}
\\

Note that the lemma proves the original conjecture of Schwarz,
namely that $M(n) \leq \binom{n}{2}$.

\begin{lemma}
\label{cycle}
Let $D$ be a critical digraph and 
$C$
a chordless cycle in it, such that 
$V(C) \neq V(D)$. 
Then
$d(v) \leq n-|V(C)|+2$ for all $v \in V(C)$, with strict 
inequality holding for at least two vertices.
\end{lemma}

{\em Proof}~
Let $J = D/V(C)$, and denote by $c$ the vertex of $J$
obtained from the contraction of $C$. Write $k$ for $|V(C)|$. The graph
$J$ has $n-k+1$ vertices, and therefore, by Lemma \ref{theorem1}, 
$d_J(c) \leq n-k$.
This implies that 
$d(v) \leq n-k+2$ for every $v \in V(C)$.

Suppose now that, for some vertex $w$ of $C$, there obtains
 $d(v) = n-k+2$ for  all vertices $v \in V(C) \setminus \{w\}$.
Then 
$d_J(c) = n-k$, 
and  all sets $E^+_D(v) \setminus V(C)$ ($v \in V(C) \setminus \{w\}$)
are equal, and the same goes for the sets $E^-_D(v) \setminus V(C)$.
Moreover, $(E^+_D(w) \setminus V(C)) \subseteq E^+_D(v)$ 
and $(E^-_D(w) \setminus V(C)) \subseteq E^-_D(v)$
for all $v \in V(C)$. But then $D - w$ must be 
strongly connected, since if $(x,w)$ and $(w,y)$ are edges
in $D$, then $y$ is reachable from $x$ in $D - w$
through vertices of $V(C) \setminus \{w\}$. {\hfill $\Box$}
\\

{\em Proof of Theorem \ref{main}}~~
The proof is by induction on $n$.
Write $s_n$ for the value claimed by Theorem \ref{main}
for $M(n)$,
namely 
$$s_n = \binom{n}{2} - n +4$$

Since $D$ is critically strongly connected, it contains a chordless cycle $C$.
 Let $|V(C)|=k$. 
If $V(C)=V(D)$ then we are done because then
$|E(D)| = n \leq s_n$. 
Thus we may assume that $V(C) \neq V(D)$,
and since $D$ is critical, 
this implies that $n \geq k+2$.

Let
$J = D/V(C)$,
and denote by $c$ the vertex of $J$
obtained from the contraction of $C$.

\begin{assertion}
\label{sn}
$$|E(D)|-|E(J)| \leq s_n-s_{n-k+1}$$
\end{assertion}

Consider first the case $k=2$.
Let $v$ be one of the two vertices of $C$.
The graph $J$, being the contraction of a strongly connected graph,
is itself strongly connected, and since $D-v$ is not strongly connected,
we have 
$J \neq D -v$.
This implies that $|E(J)| > |E(D-v)|$, and hence
$$|E(D)|-|E(J)| \leq d_D(v)-1 \leq n-2 = s_n-s_{n-1}$$

Assume now that $k \geq 3$.
Let $v$ be a vertex of $C$ having maximal degree,
namely $d_D(v) \geq d_D(u)$ for all $u \in V(C)$.
Let $r$ be the number of edges in $D$ not incident with any vertex of $C$.
Then
\[|E(D)| = r-k+\sum_{u \in V(C)} d_D(u) \]
and
\[|E(J)| = d_J(c)+r \geq d_D(v)-2+r \]
and therefore
\begin{eqnarray} 
|E(D)|-|E(J)| \leq 2-k+\sum_{u \in V(C) \setminus \{v\}} d_D(u) \nonumber
\\ \leq 2-k+2(n-k+1)+(k-3)(n-k+2) \nonumber 
\\ = (k-1)n-(k^2-2k+2) \leq s_n-s_{n-k+1} \nonumber 
\end{eqnarray}

which proves the assertion.

If $J$ is critical, then the theorem follows from Assertion \ref{sn} 
and the induction hypothesis. So, we may assume that $J$ is not critical.
But, for every vertex $u$ different from $c$, the graph $J-u$ is not 
strongly connected, since the graph $D-u$ is not strongly connected. Hence,
by Lemma \ref{theorem1}, we have

\begin{equation}
\label{djc}
d_J(c) \leq n - k
\end{equation}

On the other hand, 
the fact that $J$ is not  critical means that $J-c = D-V(C)$ is strongly 
connected.

We next show:
\begin{assertion}
\label{d_sum}
$$\sum_{v \in V(C)} d_D(v) \leq (n-1)k - n + 4 $$
\end{assertion}

{\em Proof of the assertion}~ By Lemma \ref{theorem1}
$d(v) < n$ for all vertices $v$. Hence, if $d_D(v) = 2$ for some $v \in
V(C)$ 
then 
the assertion is true. So, we may assume that 
 $d_D(v) \not = 2$ for all $v \in V(C)$. 
This means that
$(E^+_D(v) \cup E^-_D(v)) \setminus V(C) \neq \emptyset$ for every $v \in V(C)$. 
 Let $V(C) = \{v_1,v_2 \ldots v_k\}$ and $E(C) = \{(v_i, v_{i+1}): 1 \leq i
 \leq k\}$ (where, as usual, the indices are taken modulo $k$).
Without loss of generality we may assume that
$E^+_D(v_1) \setminus V(C) \neq \emptyset$. 
If $E^-_D(v_3) \setminus V(C) \neq \emptyset$ 
then $D-v_2$ is strongly connected.
Thus we may assume that $E^-_D(v_3) \setminus V(C) = \emptyset$ 
and  $E^+_D(v_3) \setminus V(C) \neq \emptyset$. 
Applying this argument again and again, we conclude that 
$k$ is even and that
$E^-_D(v_i) \setminus V(C) = \emptyset$ for all odd $i$ and
$E^+_D(v_i) \setminus V(C) = \emptyset$ for all even $i$.
By  (\ref{djc}) it follows that for every two adjacent vertices
on $C$
the total number of edges incident with them and not belonging to $C$ 
does not exceed $n-k$.
This implies that:

$$\sum_{v \in V(C)} d_D(v) \leq \frac{k}{2}(n-k)+2k \leq (n-1)k - n + 4$$

proving the assertion. 

Recall now that $n \geq k+2$ and that $D-V(C)$ is strongly 
connected. Hence $D-V(C)$ contains a chordless cycle $C'$.
Let $k' = |V(C')|$. The same arguments as above hold when $C$ is replaced by
$C'$, 
and thus we may assume that
$$\sum_{v \in V(C')} d_D(v) \leq (n-1)k' - n + 4 $$

This, together with Lemma \ref{theorem1}, yields:

$$\sum_{v \in V(D)} d_D(v) \leq (n-1)k-n+4+(n-1)k'-n+4+(n-1)(n-k-k') = 2s_n$$
which means that
$$|E(D)| \leq s_n$$

{\hfill $\Box$}

\end{document}